\documentclass[12pt]{article}
\usepackage{amsmath,amssymb}
\usepackage{amsthm}
\newtheorem{thm}{Theorem}
\newtheorem{prop}{Proposition}
\newtheorem{cor}{Corollary}
\newtheorem{lemma}{Lemma}
\newtheorem{note}{Note}
\theoremstyle{note}
\theoremstyle{remark}

\newtheorem{exam}{Example}
\theoremstyle{example}

\def\rit#1{{\mbox{\rm #1}}\hspace{1mm}}

\def\itemx#1{\item[{\rm(#1)}]}
\def\a{\mathfrak a}

\def\f{\mathfrak f}
\def\p{\mathfrak p}
\def\P{\mathfrak P}
\def\g{\mathfrak g}
\def\k{\mathfrak k}
\def\l{\mathfrak l}
\def\m{\mathfrak m}
\def\q{\mathfrak q}
\def\Q{\mathfrak Q}

\def\h{\mathfrak h}
\def\G{\mathfrak G}
\def\C{\mathfrak C}
\def\D{\mathfrak D}
\def\SS{\mathfrak S}
\def\T{\mathfrak T}
\def\rit#1{{\mbox{\rm #1}}}

\def\modx#1#2{\equiv #1 \mod ~#2}

\def\itemx#1{\item[{\rm(#1)}]}

\begin{document}
\begin{center}
{\Large{\bf Trace of Frobenius endomorphism of an elliptic curve with complex multiplication}}\footnote{Research supported by Grand-in-Aid for Scientific Research No.12640036\par
2000 {\it Mathematics Subject Clasiification}\; 11G05,14H52\par
 }\vspace{2mm}
\end{center}
\begin{center}
{\bf By  Noburo ISHII}\vspace{3mm}
\end{center}
\noindent
{\bf 1.}~Introduction\vspace{3mm}\par
Let $K=\mathbb Q(\sqrt{-m})$ be an imaginary quadratic field, where $m$ is a square-free positive integer. Let $R$ be an order of $K$ of conductor $f_0$ with a basis $\{1,\omega\}$ over $\mathbb Z$. We denote by $d(R)$ and $h(R)$ the discriminant and the class number of $R$ respectively. Let $f$ be the smallest positive integer such that $f\sqrt{-m}\in R$. Then we have $f=f_0/2$ (resp. $f_0$) if $m\modx3 4$ and $f_0$ is even (resp. otherwise). Let $E$ be an elliptic curve with complex multiplication by $R$ and denote by $j(E)$ the $j$-invariant of $E$. We may assume that $E$ is defined by a Weierstrass equation:~$y^2=x^3+Ax+B,A,B\in F=\mathbb Q(j(E))$. First, we introduce the notation used in the follwing.  For an endomorphism $\lambda$ of $E$, the kernel of $\lambda$ is denoted by $E[\lambda]$. For a prime ideal $\p$ of $F$, we denote by ${\ell_\p}$ the relative degree of $\p$ over $\mathbb Q$. If $E$ has good reduction at $\p$, then we denote by $\overset{\sim}{E_\p}$ the reduction of $E$ modulo $\p$. For a point $P$ of $E$ we denote by $P^\sim$ the reduction of $P$ modulo $\p$. Further we denote by $\varphi_\p$ the Frobenius endomorphism of $\overset{\sim}{E_\p}$ and by $a_\p(E)$ the trace of $\varphi_\p$. By $\mathbb F_q$, we denote the finite field of $q$-elements.  If $\overset{\sim}{E_\p}$ is defined over $\mathbb F_q$, then $\overset{\sim}{E_\p}(\mathbb F_{q})$ denotes the group of $\mathbb F_{q}$-rational points of $\overset{\sim}{E_\p}$.\par
 Now let $p$ be an odd prime number and $\p$ a prime ideal of $F$ dividing $p$. Let us assume that $p$ and $\p$ satisfy the following condition:\vspace{2mm}\newline
(1) {\it $p$ splits completely in $K$, $p$ is prime to $f$ and $E$ has good reduction at $\p$.}\vspace{2mm}\newline
 Then by complex multiplication theory (cf.~II of Silverman [11]), we know that $E$ has ordinary good reduction at $\p$ and the endomorphism ring of $\overset{\sim}{E_\p}$ is isomorphic to $R$ (cf. Theorem 12 of 13.4 of Lang [6]). Further $K(j(E))$ is the ring class field of $K$ of conductor $f_0$  (cf.~\S 9 of Cox [3]). Since $\p$ is of relative degree ${\ell_\p}$, there exist positive integers $u_p$ and $v_p$ such that 
$$4p^{{\ell_\p}}=u_p^2+mf^2v_p^2,~(u_p+v_p f\sqrt{-m})/2 \in R,~(u_p,p)=1. $$ 
 By the assumption, we may write $\varphi_{\p}=(a_\p(E)+b_{\p}(E)f\sqrt{-m})/2=\alpha+\beta\omega$, where $b_\p(E), \alpha$ and $\beta$ are integers. It is known that the group $\overset{\sim}{E_\p}(\mathbb F_{p^{\ell_\p}})$ is of order $N_\p(E)=p^{\ell_\p} +1-a_\p(E)$ and is isomorphic to the group $\mathbb Z/(N_\p(E)/d)\mathbb Z\oplus \mathbb Z/d\mathbb Z$, where $d$ is the greatest common divisor of $\alpha-1$ and $\beta$. On the other hand, if $d(R)<-4$, then we have easily $a_\p(E)=\epsilon_{\p}u_p$, where $\epsilon_{\p}=1\text{ or }-1$. It is easy to find $u_p$ for a given number $4p^{{\ell_\p}}$ such that $4p^{{\ell_\p}}=u_p^2+mf^2v_p^2,~(u_p,p)=1$. Therefore if we determine $\epsilon_{\p}$, then we can compute the numbers $N_\p(E)$ and $d$ rapidly. The problem to detemine $\epsilon_{\p}$ in the case $h(R)=1$ has been solved by A.R. Rajwade, A. Joux-F. Morain and others. See A. Joux and F. Morain[5] for the references of their results. In the case $h(R)=2$, this problem is solved only for one case of the order of discriminant $-20$, by  F. Lepr\'{e}vost and F. Morain [7], using the results of B.W. Brewer [1,2] for the character sum of Dickson polynomial of degree $5$. \par 
The purpose of this article is to determine $\epsilon_{\p}$ for an elliptic curves $E$ having complex multiplication by $R$ and for prime ideals $\p$ of $F$ satisftying (1), where $R$ are orders such that $h(R)=2\text{ or }3$ and $mf^2$ is divided by at least one of $3,4$ and $5$.  Thus $R$ are orders of discriminant 
\begin{align*}
d(R)=&-15,-20,-24,-32,-35,-36,-40,-48, -51,-60,-64,-72,-75,\\
&-99,-100,-108,-112,-115,-123,-147,-235,-243, -267. 
\end{align*}
Further we assume that $j(E)$ is real to avoid tedious argument.  \par
 Our idea to solve the problem is as follows (for details see \S 2). Let $s$ be a divisor of $f^2m$ and assume $s\geq 3$. We find a $F$-rational cyclic subgroup $C_s$ of $E[f\sqrt{-m}]$ of order $s$  and take a generator $Q$ of $C_s$.   Consider the Frobenius isomorphism of $\sigma_\p$ of ${\p}$. Then $F$-rationality of $C_s$ shows $Q^{\sigma_\p}=[r_{\p}](Q)$ for an integer $r_\p$.  Using $Q^\sim\in\overset{\sim}{E_\p}[f\sqrt{-m}]$ and $(Q^{\sigma_\p})^\sim =\varphi_\p(Q^\sim)$, we have
 \begin{align*}
 [2]([r_{\p}](Q))^\sim &=[2](Q^{\sigma_\p})^\sim =[2]\varphi_\p(Q^\sim) \\
&=[(a_\p(E)+b_{\p}(E)f\sqrt{-m})](Q^\sim)= [a_\p(E)](Q^\sim).
\end{align*}
This shows that $a_\p(E)\modx {2r_{\p}}s$. Therefore the number $\epsilon_{\p}$ is determined by the condition $\epsilon_{\p}u_p\modx {2r_{\p}}s$. This argument reduces our original problem to a problem of finding a point $Q$ and of determining $r_\p$ for a given prime ideal $\p$. In \S 2, we give auxiliary results to find the cyclic sbgroup $C_s$ and a generator $Q$. If $s$ is an odd prime number, then we show, in Proposition~\ref{prop3} of \S 2, that the $s$-division polynomial $\Psi_s(x,E)$ of $E$ has a unique $F$-rational factor $H_{1,E}(x)$ of degree $(s-1)/2$ and that the point $Q$ is obtainable from a solution of $H_{1,E}(x)=0$. In \S 3 we determine $r_\p$ for the case $f^2m$ is divided by $3\text{ or }4$ and in \S 4 for the case $f^2m$ is divided by $5$. Though we deal with a specified elliptic curve $E$ for each order $R$, a similar result is easily obtained for any elliptic curve $E'$ of the $j$-invariant $j(E)$, because $E'$ is a quadratic twist of $E$ and $a_\p(E')$ is the product of $a_\p(E)$ and the value at $\p$ of the character corresponding to the twist. \par 
 In the following, we assume any elliptic curve is defined by a short Weierstrass equation.
\vspace{3mm}\newline 
\noindent
{\bf 2.}~The subgroups of $E[f\sqrt{-m}]$ and decomposition of division polynomials\vspace{3mm}\par\noindent
{\bf 2.1.}~Let $E$ be an elliptic curve with complex multiplication by $R$. By the definition of $f$, we have $f\sqrt{-m} \in R$. 
\begin{prop}\label{prop1} The group $E[f\sqrt{-m}]$ is  cyclic of order $f^2m$.
\end{prop}
\begin{proof} By Proposition 2.1 of Lenstra [8], we know $E[f\sqrt{-m}]$ is isomorphic to $R/f\sqrt{-m}R$. Let $f$ be odd and $m\modx 34$. Then $R=\mathbb Z\oplus f\omega\mathbb Z$, where $\omega =(1+\sqrt{-m})/2$. Further $f\sqrt{-m}R=f(2\omega -1)\mathbb Z\oplus f^2(\omega -(m+1)/2)\mathbb Z$. Put $\xi =f\omega -f(mf+1)/2\in f\sqrt{-m}R$. Then we have
$f(2\omega -1)=2\xi+mf^2,~f^2(\omega-(m+1)/2)=f\xi+(f-1)f^2m/2$. This shows that $\{f^2m,\xi\}$ is a basis of $f\sqrt{-m}R$ over $\mathbb Z$. Since $\{1,\xi\}$ is a basis of $R$ over $\mathbb Z$, $R/f\sqrt{-m}R$ is a cyclic group of order $f^2m$. The assertion for the other cases is easily obtained. 
 \end{proof}
 \begin{lemma}\label{lem1} Let $r$ be a fixed prime number. Then there exist infinitely many prime numbers of the form $u^2+v^2f^2m$, where $u$ and $v$ are integers and $v$ is prime to $r$.
\end{lemma}
\begin{proof} Consider the ideal groups $G_0$ and $P_0$ of $K$ such that
\begin{equation*}
G_0=\{{\mathfrak a}\mid {\mathfrak a}\text{ is prime to }2r fm\},~
P_0=\{(\alpha )\mid \alpha \modx 1{2rf\sqrt{-m}}\}.
\end{equation*}
Then $P_0$ is a subgroup of $G_0$ of finite index and by Tshebotareff's density theorm, in each factor class there exist infinitely many prime ideals of degree $1$. Let $\gamma=u_0+v_0f\sqrt{-m}$ such that ideal $(\gamma) \in G_0$ and $u_0,v_0\in \mathbb Z$ and further  $v_0$ is prime to $r$. Then every integral ideal of the class $(\gamma)P_0$ has a generator of the form $u_1+v_1f\sqrt{-m}~(\,u_1,v_1\in\mathbb Z,r\nmid v_1)$. Thus 
we have our assertion.
\end{proof}
In the following, let $p$ be an odd prime number and $\p$ a prime ideal of $F$ dividing $p$ and assume that $p$ and $\p$ satisfy the condition (1).
\begin{lemma}\label{lem2} Let $s$ be an odd prime number dividing $f^2m$. Let $q=p^{\ell_\p}$. Assume that  $q=u^2+v^2f^2m,~(v,ps)=1$ or $4q=u^2+v^2f^2m,~(v,2sp)=1$.  Then we have
$$\overset{\sim}{E_\p}[s]\cap \overset{\sim}{E_\p}[f\sqrt{-m}]\setminus \{0\}=\{P=(\alpha,\beta)\in \overset{\sim}{E_\p}[s]~\arrowvert~s\nmid [\mathbb F_q(\alpha):\mathbb F_q ]\},$$
where $[\mathbb F_q(\alpha):\mathbb F_q ]$ denotes the degree of the field $\mathbb F_q(\alpha)$ over $\mathbb F_q$.
\end{lemma}
\begin{proof}
By the assumption, $\overset{\sim}{E_\p}$ is defined over $\mathbb F_q$. First we assume the Frobenius endomorphism $\varphi_\p$ is given by $\varphi_\p=(u+vf\sqrt{-m})/2$, if necessary, after replacing $u$ by $-u$ or $v$ by $-v$. Let $P=(\alpha,\beta)\in \overset{\sim}{E_\p}[s]$. If $P\in \overset{\sim}{E_\p}[f\sqrt{-m}]$, then, for $h=(s-1)/2$, we have $\varphi^h_\p([2^h](P))=[u^h](P)$. Since $2^h,u^h\modx{\pm 1}s$, we have $\varphi^h_\p(P)=\pm P$. This shows $[\mathbb F_q(\alpha):\mathbb F_q ]\le (s-1)/2$. Conversly let $P=(\alpha,\beta)\in \overset{\sim}{E_\p}[s], s\nmid k=[\mathbb F_q(\alpha):\mathbb F_q]$ and $r=q^k$. Since $\varphi^k_\p(P)=(\alpha^r,\beta^r)=(\alpha,\beta^r)=\epsilon P~(\epsilon =\pm 1)$,we have $[(u+vf\sqrt{-m})/2)^k-\epsilon](P)=0$. Since $P\in \overset{\sim}{E_\p}[s]$ and $s|f^2m$, we have $[(u^k-2^k\epsilon)+ku^{k-1}vf\sqrt{-m}](P)=0$ and $[(u^k-2^k\epsilon)^2+(ku^{k-1}v)^2f^2m](P)=0$. Thus $[(u^k-2^k\epsilon)^2](P)=0$. Since the order of $P$ is $s$, we see $[(u^k-2^k\epsilon)](P)=0$ and $[ku^{k-1}vf\sqrt{-m}](P)=0$. By the assumption, $k,u$ and $v$ are prime to $s$. Therefore we conclude $[f\sqrt{-m}](P)=0$. Hence $P\in \overset{\sim}{E_\p}[f\sqrt{-m}]$. In the case $\varphi_\p=u+vf\sqrt{-m}$, the same argument holds true.
\end{proof}
\begin{cor}\label{cor1}  Let $\Psi_s(x,\overset{\sim}{E_\p})$ be the $s$-division polynomial of $\overset{\sim}{E_\p}$. Then we know $\Psi_s(x,\overset{\sim}{E_\p})$ is the product of two $\mathbb F_q$-rational polynomials $h_1(x)$ and $h_{2}(x)$ such that $h_{1}(x)$ is of degree $(s-1)/2$ and the degree of every irreducible factor of $h_{2}(x)$ is divided by $s$. Further the solutions of $h_{1}(x)=0$ consist of all distinct $x$-coordinates of non-zero points in $\overset{\sim}{E_\p}[s]\cap \overset{\sim}{E_\p}[f\sqrt{-m}]$.
\end{cor}
\begin{proof} Since $p$ is prime to $f^2m$, by Proposition~\ref{prop1}, $\overset{\sim}{E_\p}[s]\cap \overset{\sim}{E_\p}[f\sqrt{-m}]$ is a $\mathbb F_q$-rational cyclic group of order $s$. Thus if we put $h_{1}(x)=\prod_\alpha(x-\alpha)$, where $\alpha$ runs over all distinct $x$-coordinates of non-zero points in $\overset{\sim}{E_\p}[s]\cap \overset{\sim}{E_\p}[f\sqrt{-m}]$, then $h_{1}(x)$ is $\mathbb F_q$-rational and of degree $(s-1)/2$. The assertion for $h_{2}(x)$ follows immediately from Lemma~\ref{lem2}.\end{proof}
\begin{lemma}\label{lem3} Let $4|f^2m$ and $q=p^{\ell_\p}=u^2+v^2f^2m,~(v,2)=1$. Let $Q_1$ be an point of order $4$ of $\overset{\sim}{E_\p}[f\sqrt{-m}]$ and $Q_2$ a point of $\overset{\sim}{E_\p}$ such that $[2](Q_1)=[2](Q_2)$ and $Q_2\ne\pm Q_1$. Then the $x$-coordinates $x_1$ and $x_2$ of $Q_1$ and $Q_2$ are all $\mathbb F_q$-rational solutions of $\Psi_4(x,\overset{\sim}{E_\p})/y=0$. Furthermore let $y^2=h(x)$ be the equation of $\overset{\sim}{E_\p}$. Assume that $\varphi_\p=u+vf\sqrt{-m}$. Then, of two elements $x_1$ and $x_2$, only $x_1$ satisfies the relation $(h(x_1)/\p)=(-1)^{(u-1)/2}$,where $(\phantom{a}/\p)$ denotes the Legendre symbol for $\p$.
\end{lemma}
\begin{proof} Since $\overset{\sim}{E_\p}[f\sqrt{-m}]$ is $\mathbb F_q$-rational, we see $x_1$ and $x_2$ are $\mathbb F_q$ -rational. Let $\alpha$ be a $\mathbb F_q$-rational root of $\Psi_4(x,\overset{\sim}{E_\p})/y=0$ and put $S=(\alpha,\beta)$. Then $S$ is a $4$-division point of $\overset{\sim}{E_\p}$ and we have $$\varphi_\p(S)=[u+vf\sqrt{-m}](S)=(\alpha^q,\beta^q)=(\alpha,\pm\beta)=[\varepsilon](S),\medskip (\varepsilon =\pm 1).$$
Thus we have $[(u-\varepsilon)+vf\sqrt{-m}](S)=0$. This shows $[(u-\varepsilon)^2+v^2f^2m](S)=0$. Since the order of $S$ is $4$, $u-\varepsilon$ is divided by $2$. Thus $[f\sqrt{-m}]([2]S)=0$. Since $[2](Q_1)$ is the only one point of degree $2$ in $\overset{\sim}{E_\p}[f\sqrt{- m}]$, we have $[2](S)=[2](Q_1)$. This shows that $S$ equals to one of $\pm Q_1$ and $\pm Q_2$. Therefore $\alpha$ equals to $x_1$ or $x_2$. Let $P=(x,y)$ be a point of $\overset{\sim}{E_\p}$ of order $4$ such that $x \in\mathbb F_q$. Then
\begin{align*}
\varphi_\p(P)&=(x^q,y^q)=(x,yh(x)^{(q-1)/2})\\
&=[(h(x)/\p)](P)=[u](P)+[vf\sqrt{-m}](P).
\end{align*}
Therefore we have $(h(x)/\p)\modx u 4$ if and only if $P\in \overset{\sim}{E_\p}[f\sqrt{- m}]$.
\end{proof}
{\bf 2.2.}~Let $s$ be a positive divisor of $f^2m$ and $s\ge 3$. By Proposition~\ref{prop1}, there exists a unique subgroup $C_s$ of $E[f\sqrt{-m}]$ of order $s$. Let $Q=(x_Q,y_Q)$ be a generator of $C_s$. Consider the fields $L=F(x_Q)$ and $M=F(Q)$. Since $E[f\sqrt{-m}]$ is $F$-rational, $C_s$ is $F$-rational and the field $M$ is an abelian extension over $F$. By class field theory, the Galois group $G$ of $M$ over $F$ is isomorphic to an ideal class group $\G$ of $F$. For an ideal class $\C\in\G$, let $\sigma_\C$ be the isomorphism of $\G$ corresponding to $\C$.  Then we have
\begin{thm}\label{th1}  Let $\C$ be the class represented by $\p$ and $Q^{\sigma_\C}=[r_\C](Q)$. Then we have $a_\p(E)\modx{2r_\C}s$. Further if $a_\p(E)$ is even, then we have $a_\p(E)/2\modx{r_\C}s$. 
\end{thm}
\begin{proof} Let $\varphi_\p=(a_\p(E)+b_\p(E) f\sqrt{-m})/2$. Since $(Q^{\sigma_\C})^\sim=\varphi_\p(Q^\sim)$, we see
\begin{equation*}
\begin{split}
[2]([r_\C](Q))^\sim &=[2](Q^{\sigma_\C})^\sim =[2]\varphi_\p(Q^\sim) \\
&=[(a_\p(E)+b_\p(E) f\sqrt{-m})](Q^\sim)= [a_\p(E)](Q^\sim).
\end{split}
\end{equation*}
Since $p$ is prime to $s$, $Q^\sim$ is of order $s$. Thus $a_\p(E)\modx{2r_\C}s$. If $a_\p(E)$ is even, then $\varphi_\p=(a_\p(E)/2)+(b_\p(E)/2) f\sqrt{-m}$. By a similar argument we have $[a_\p(E)/2](Q^\sim)=[r_\C](Q)$. This shows the remaining assertion.
\end{proof}
\begin{prop}\label{prop2} Let $s$ be an odd prime number of $f^2m$. If $p^{\ell_\p}\modx 1s$, then $$a_\p(E)\modx{2(y_Q^2/\p)}s$$.
\end{prop}
\begin{proof} Since we have $4p^{\ell_\p}=a_\p(E)^2+b_\p(E)^2f^2m$, Theorem~\ref{th1} shows that $r_\C\modx {\pm 1}s$. Thus $x_Q^\sim\in\mathbb F_q$. By the similar argument in the last part of Lemma~\ref{lem3}, we have our assertion.\end{proof}
\begin{prop}\label{prop3} Let $s$ be an odd prime divisor of $f^2m$ and $\Psi_s(x,E)$ the $s$-division polynomial of $E$. Then $\Psi_s(x,E)$ is the product of two $F$-rational polynomials $H_{1,E}(x)$ and $H_{2,E}(x)$ such that $H_{1,E}(x)$ is of degree $(s-1)/2$ and every irreducible factor of $H_{2,E}(x)$ is of degree $s$. Further the solutions of $H_{1,E}(x)=0$ consist of all distinct $x$-coordinates of non-zero points of $C_s$.
\end{prop}
\begin{proof} Let $H_{1,E}(x)=\displaystyle\prod_{t}(x-t)$, where $t$ runs over all distinct  $x$-coordinates of non-zero points of $C_s$. Since $C_s$ is $F$-rational, we see $H_{1,E}(x)$ is $F$-rational of order $(s-1)/2$ and clearly it divides $\Psi_s(x,E)$. By Lemma~\ref{lem1}, we can choose an odd prime $p$ and a prime ideal $\p$ dividing $p$ such that they satisfy (1) and $p$ is of the form $p=u^2+v^2f^2m,~(v,s)=1$, and further the reduction of $\Psi_s(x,E)$ modulo $\p$ is equal to $\Psi_s(x,\overset{\sim}{E_\p})$. Take a point $P\in E[s]\setminus C_s$ and put $Q=[f\sqrt{-m}](P)$.  Clearly, we have $Q\in E[f\sqrt{-m}]$ and $E[s]=<P>\oplus <Q>$. Let $G_1$ be the Galois group of $F(E[s])$ over $F$. By the representation of $G_1$ on $E[s]$ with the basis $\{P,Q\}$, $G_1$ is identified with a subgroup of the group
$$G_0=\left\{\left.\begin{pmatrix}a & 0 \\ b & \pm a\end{pmatrix}\right| a\in\mathbb F^\times_s,b\in\mathbb F_s\right\}.$$  Consider the subgroups of $G_0$:
$$H=\left\{\left.\begin{pmatrix}a & 0 \\ 0 & \pm a\end{pmatrix}\right| a\in\mathbb F^\times_s\right\},~N=\left\{\left.\begin{pmatrix}1 & 0 \\ b & 1\end{pmatrix}\right|b\in\mathbb F_s\right\}.$$
Then we see $G_0=HN,~G_0\triangleright N$ and $H\cap N=\{1_2\}$, where $1_2$ is the unit matrix. Since the order of $N$ is $s$ and $s$ is prime, we know that $G_1\supset N$ or $G_1\cap N=\{1_2\}$. Let $\Omega$ be the set of all subgroups of order $s$ of $E[s]$. Then $\Omega$ consists of $s+1$ elements and $G_1$ operates on $\Omega$. By Corollary 1, the degree of every irreducible factor of $H_{2,E}(x)=\Psi_s(x,E)/H_{1,E}(x)$ is divided by $s$. Thus we know $C_s$ is one and the only one fixed point of $G_1$.  First let us consider the case $G_1\supset N$. Then we have $G_1=H_1N,~H_1=H\cap G_1$. Since $H_1$ is the fixed subgroup of $<P>$, the orbit of $<P>$ consists of $s$ elements. Therefore $\Omega$ decomposes into two orbits. In particular, for each $n,~1\le n \le (s-1)/2$, the $x$-coordinate of $[n]P$ has $s$ conjugates over $F$. Thus every irreducible factor of $H_{2,E}(x)$ is of degree $s$. Next consider the case $G_1\cap N=\{1_2\}$. Then the order of $G_1$ is a divisor of $2(s-1)$ and is prime to $s$. Since the order of a matrix $\begin{pmatrix}a&0\\b&a\end{pmatrix},~(b\ne 0)$ is divided by $s$, $G_1$ dose not contain the matrices of this form. Therefore there exists a $\lambda\in\mathbb F_s$ such that $G_1$ is contained in the subgroup 
$$\left<\alpha,\left.\begin{pmatrix}1&0\\ \lambda & -1\end{pmatrix}~\right |~\alpha\in\mathbb F_s^\times\right >.$$
This shows $<P+(\lambda /2)Q>$ is a fixed point. Thus we have a contradiction.
\end{proof}
\begin{prop}\label{prop4} Let $4|f^2m$. If $Q$ is a point of order $4$ in $E[f\sqrt{-m}]$ and $T$ is a poinf of $E$ such that $[2](Q)=[2](T)$ and $T\ne\pm Q$, then the $x$-coordinates $x_{Q}$ and $x_{T}$ of $Q$ and $T$ are all $F$-rational solutions of $\Psi_4(x,E)/y=0$.
\end{prop}
\begin{proof} Using Lemma~\ref{lem3} instead of Lemma~\ref{lem2} and tracing the argument in the first part of Proposition~\ref{prop2}, we have the assertion.
\end{proof}
 In \S 4, to study the ideal class groups of $F$ corresponding to the fields $L$ and $M$, we must determine conductors $\f_L$ and $\f_M$ of abelian extensions $L$ and $M$ over $F$.  In next lemma, we shall give some results for the conductors. For a prime ideal $\q$ and an integral ideal $\a $ of $F$, we denote by $e_\q(\a)$ the maximal integer $m$ such that $m\ge 0$ and $\q^m$ dividing $\a$. 
\begin{lemma}\label{lem4} Let $Q$ be a point of $E$ of order $s$.  Assume that $s$ is an odd prime number, $s>3$ and $Q$ generates a $F$-rational subgroup $<Q>$.  Let $L,M,\f_L$ and $\f_M$ be as above. If $\q$ is a prime ideal of $F$ prime to $(2s)$, then $e_\q(\f_L)\le e_\q(\f_M)$ and $e_\q(\f_M)>0$ implies $e_\q(\f_L)>0$. Further if $E$ has good reduction at $\q$, then $e_\q(\f_M)=0$.
\end{lemma} 
\begin{proof} Since $L$ is a subfield of $M$, we have $e_\q(\f_L)\le e_\q(\f_M)$. If $E$ has good reduction at $\q$, then N\'{e}ron-Ogg-Shafarevich criterion (cf.Proposition 4.1 of VII of Silverman [10]) shows that $e_\q(\f_M)=0$. We shall prove $e_\q(\f_M)>0$ implies $e_\q(\f_L)>0$. Assume that $e_\q(\f_M)>0$ and $e_\q(\f_L)=0$, thus, assume that $\q$ is ramified in $M$ and is unramified in $L$. Let $\Q$ be a prime ideal of $M$ dividing $\q$ and $M_\Q$ the completion of $M$ with respect to $\Q$.  Further we denote by $k_M$ the residue field of $\Q$. Let $E_0,E_1$ and $\overset{\sim}{E}_{ns}$ be the groups defined in VII of Silverman [10]. Since $E$ has additive reduction at $\q$, by Theorem 6.1 of VII of [10] we have $[E(M_\Q):E_0(M_\Q)]=w\le 4$. Since $Q$ has order $s$, by replacing $Q$ by $[w] Q$ if necessary, we can assume that $Q\in E_0(M_\Q)$. Let $\sigma$ be a non trivial element of inertia group of $\Q$. Then since $x_Q^\sigma=x_Q$, we have $Q^\sigma=-Q$. By considering the reduction modulo $\q$, we have $Q^{\sim}=-Q^{\sim}$. Therefore $Q^\sim\in \overset{\sim}{E}_{ns}(k_M)$ and $[2](Q^\sim)=0$. Since the characteristic of $k_M$ is prime to $2$, Proposition 5.1 of VII of [10] shows $Q^\sim=0$. Therefore by Proposition 2.1 of VII of [10] , we know $Q\in E_1(M_\Q)$. Consequently, by Proposition 3.1 of VII of [10], we have $Q=0$. This contradicts that $Q\ne 0$. 
\end{proof}

Finally for $s=3,4,5$, we list $s$-division polynomials $\Psi_s(x,E)$:
\begin{equation*}
\begin{cases}
\Psi_3(x,E)&=3x^4+6Ax^2+12Bx-A^2,\\
\Psi_4(x,E)&=2y(2x^6+10Ax^4+40Bx^3-10A^2x^2-8ABx-16B^2-2A^3),\\
\Psi_5(x,E)&=5x^{12}+62Ax^{10}+380Bx^9-105A^2x^8+240ABx^7\\
&-(300A^3+240B^2)x^6-696A^2Bx^5-(125A^4+1920AB^2)x^4\\
&-(1600B^3+80BA^3)x^3-(50A^5+240A^2B^2)x^2\\
&-(640AB^3+100A^4B)x+A^6-32B^2A^3-256B^4.
\end{cases}
\end{equation*}\vspace*{3mm}
\noindent
{\bf 3.}~ The case $f^2m$ is divided by $3$ or $4$\par
Let $s=3$ or $4$. Assume that $s|f^2m$. Let $Q=(x_Q,y_Q)$ be a point of $E[f\sqrt{-m}]$ of order $s$ . By Propositions~\ref{prop2} and~\ref{prop3}, we know $x_Q\in F$. We may write $y_Q^2=w^2\alpha_E$ such that $w \in F^\times$, $\alpha_E$ is an integer of $F$ and ideal$(\alpha)$ has no square factors. In the following, let $p$ be an odd prime number and $\p$ a prime ideal of $F$ dividing $p$ and assume they satisfy the condition (1). Then there exist positive integers $u_p$ and $v_p$ such that $4p^{\ell_\p}=u_p^2+mf^2v_p^2,~(u_p+v_pf\sqrt{-m})/2\in R, ~(u_p,p)=1$. If $u_p$ is even, then clearly we have $p^{\ell_\p}=(u_p/2)^2+mf^2(v_p/2)^2$. 
\begin{thm}\label{th2}  Let $u_p$ and $v_p$ be as above. If we choose $\epsilon_\p\in\{\pm 1\}$ such that 
$\epsilon_\p (u_p/2)\modx {(\alpha_E/\p)} s$, then we have $a_\p(E)=\epsilon_\p u_p$. \end{thm}
\begin{proof} Since $F(Q)=F(\sqrt{\alpha_E})$, we have $Q^{\sigma_\p} =[(\alpha_E/\p)](Q)$. Thus by Theorem~\ref{th1}, we have our assertion. It is noted that $u_p$ can be odd only in the case $s=3$.
\end{proof}
Let $E_0$ be an elliptic curve defined by a Weierstrass equation:~$y^2=x^3+A_0x+B_0~(A_0,B_0\in F$). If $E_0$ is isomorphic to $E$ over an extension $F_0$ over $F$, then there exists an element $\delta \in F_0$ such that $A_0=\delta ^4A,~B_0=\delta^6B$. Since $j(E)\ne 0,1728$, we know that $A,B,A_0$ and $B_0$ are not $0$ and $\delta^2\in F$. Therefore we may put $\alpha_{E_0}=\delta^2\alpha_E$. In particular we obtain
\begin{thm}\label{th3}  Let $E^*$ be the twist of $E$ defined by the equation $y^2=x^3+A\alpha_E^2+B\alpha_E^3$. Further assume that $E^*$ has good reduction at $\p$. Let $u_p$ and $v_p$ be as above. If we choose $\epsilon_\p\in\{\pm 1\}$ such that $\epsilon_p (u_p/2)\modx 1 s$, then we have $a_\p(E^*)=\epsilon_\p u_p$.  \end{thm}
\noindent 

 The $j$-invariants of elliptic curves with complex multiplication by $R$ are solutions of the class equation  $H_{|d(R)|}(x)=0$ of discrimint $d(R)$ (cf. \S 13 of Cox [3] for the class equations). In the following, we shall use the table of class equations prepared by M.Kaneko. We shall gives a canonical elliptic curve $E$ with complex multiplication by $R$ and compute $\alpha_E$ in subsections 3.1 and 3.2 for the cases $s=3$ and $4$ respectively.   \vspace{3mm}\newline 
{\bf 3.1.} The case $s=3$\vspace{3mm}\par
 We shall explain the process to obtain a canonical elliptic curve $E$ in the case $d(R)=-15$. At first we take a solution $j_1=(-191025+85995\sqrt{5})/2$ of the equation:
$$H_{15}(x)=x^2+191025x-121287375=0.$$ 
Let $E_1$ be the elliptic curve defined by the equation: 
$$y^2=x^3+A_1x+B_1,\, A_1=-1/48-36/(j_1-1728), B_1=1/864+2/(j_1-1728).$$
Then the $j$-invariant of $E_1$ is equal to $j_1$.
 By considering twists of $E_1$ by elements $\sqrt n,~n\in F=\mathbb Q(\sqrt 5)$, we find an elliptic curve $E$ such that coefficients $A$ and $B$ of an equation~$y^2=x^3+Ax+B$ of $E$ are integers of $F$ and further the absolute value of the norm of the square factor of $A$ is as small as possible. In this case, we take $n=2^237(4+\sqrt{5})/(\sqrt{5}(4-\sqrt{5}))$. Therefore we see $A=A_1n^2=105+48\sqrt{5},B=B_1n^3=-784 -350\sqrt{5}$ and  
$$\Psi_3(E,x)=3(x^3+6x^2+3\sqrt{5}x^2+(291+132\sqrt{5})x+590+265\sqrt{5})(x-6-3\sqrt{5}).$$ 
This shows $x_Q=6+3\sqrt{5}$ and $y_Q^2=2^4((1+\sqrt{5})/2)^{11}$. Finally we have
\begin{prop}\label{prop5} Let $E$ be the elliptic curve defined by the equation 
$$y^2=x^3+(105+48\sqrt{5})x-784 -350\sqrt{5}.$$
 Then $E$ has complex multiplication by the order of discriminant $-15$. Further we have $\alpha_E=(1+\sqrt{5})/2$.
\end{prop}
\begin{note}
For another root $j_2=(-191025-85995\sqrt{5})/2$ of $H_{15}(x)=0$, we consider the conjugate elliptic curve $\overline E$ of $E$ over $\mathbb Q$ and put $\alpha_{\overline E}= (1-\sqrt{5})/2$.\end{note}
\vspace{3mm}
\noindent
 Example 1.\newline
(1)~ Let $p=61$. Then $(-15/p)=(5/p)=1$. Thus ${\ell_\p} =1$. Choose the  prime ideal $\p$ such that $\p\ni\sqrt 5-26$. Since $(\alpha_E/\p)=(54/61)=-1$ and $4p=2^2+4^215$, $a_\p(E)=-2.$\newline
(2)~ Let $p=83$. Then $(-15/p)=1,(5/p)=-1$. Thus ${\ell_\p} =2$. Since $(\alpha_E/(p))=-1$ and $4p^2=154^2+16^2\cdot 15$, $a_{(p)}=154$.\vspace{3mm}\par
 For other cases, we give only results and deta nacessary to obtain the results.  For each order $R$, the data consits of the class polynomial $H_{|d(R)|}(x)$, a solution $j$ of $H_{|d(R)|}=0$, coefficients $A$ and $B$ of a Weierstrass equation of an elliptic curve $E$ with $j(E)=j$, $x_Q$, $y_Q^2$ and $\alpha_E$. We list them in the follwing format:
 \begin{center}
\begin{tabular}{|c|c|} \hline
$d(R)$ & $\phantom{aaa}H_{|d(R)|}(x)\phantom{aaa}$\\\cline{2-2}
           & $j$\\\cline{2-2}
           & $A,~ B$\\\cline{2-2} 
           &$x_Q,~y_Q^2$\\\cline{2-2}
           &$\alpha_E$\\\hline
\end{tabular}
\end{center}
\begin{center}
\begin{tabular}{|c|c|}
\multicolumn{2}{l}{The results and data for the case $h(R)=2$}\\\hline

$-24$ & $x^2-4834944x+14670139392$\\\cline{2-2}
           & $2417472+1707264\sqrt 2$\\\cline{2-2}
           & $-21+12\sqrt 2,~ -28+22\sqrt 2$\\\cline{2-2} 
           &$-3+3\sqrt{2}, ~2(1-\sqrt 2)^6(1+\sqrt 2)$\\\cline{2-2}
           &$1+\sqrt 2$\\\hline
$-36$ & $x^2-153542016x-1790957481984 $\\\cline{2-2} 
      &$76771008+44330496\sqrt 3 $\\\cline{2-2}
      &$-120-42\sqrt{3}, ~448 +336\sqrt{3}$\\\cline{2-2}
      &$3+3\sqrt{3}, ~4(2+\sqrt{3})^2(1+\sqrt{3})$\\\cline{2-2}
      &$1+\sqrt{3}$\\\hline 
      $-48  $ & $ x^2-2835810000x+6549518250000    $\\\cline{2-2}
           & $1417905000+818626500\sqrt{3}   $\\\cline{2-2}
           & $-1035-240\sqrt{3}, ~12122 +5280\sqrt{3}   $\\\cline{2-2} 
           &$-9+18\sqrt{3}, ~4(2-\sqrt{3})^4(1-2\sqrt{3})^2(8+6\sqrt{3})$\\\cline{2-2}
           &$8+6\sqrt{3} $\\\hline 
$-51  $ & $x^2+5541101568x+6262062317568     $\\\cline{2-2}
           & $-2770550784-671956992\sqrt{17}   $\\\cline{2-2}
           & $-60-12\sqrt{17}, ~-210-56\sqrt{17}   $\\\cline{2-2} 
           &$ -6,~ -2(4-\sqrt{17})^2 $\\\cline{2-2}
           &$ -2$\\\hline            
$-60  $ & $x^2-37018076625x+153173312762625     $\\\cline{2-2}
           & $(37018076625+16554983445\sqrt{5})/2   $\\\cline{2-2}
           & $(-645+201\sqrt{5})/2, ~1694-924\sqrt{5}   $\\\cline{2-2} 
           &$-(45-15\sqrt{5})/2, ~-((1-\sqrt{5})/2)^{16}  $\\\cline{2-2}
           &$-1 $\\\hline 
$-72  $ & $ x^2-377674768000x+232381513792000000    $\\\cline{2-2}
           & $188837384000+77092288000\sqrt 6,   $\\\cline{2-2}
           & $-470-360\sqrt 6, ~19208+10080\sqrt{6}   $\\\cline{2-2} 
           &$6+9\sqrt 6, ~4(5-2\sqrt 6)^2(2+\sqrt 6)  $\\\cline{2-2}
           &$2+\sqrt 6 $\\\hline 
 $-75  $ & $x^2+654403829760x+5209253090426880     $\\\cline{2-2}
           & $-327201914880+146329141248\sqrt{5}   $\\\cline{2-2}
           & $-2160+408\sqrt{5},~42130-10472\sqrt{5}   $\\\cline{2-2} 
           &$-(15+21\sqrt{5}),~ (-25-13\sqrt{5})(4-\sqrt{5})^2((1+\sqrt{5})/2))^{14}  $\\\cline{2-2}
           &$-25-13\sqrt{5} $\\\hline            
\end{tabular}          
\begin{tabular}{|c|c|}\hline
 $-99  $ & $x^2+37616060956672x-56171326053810176  $\\\cline{2-2}
           & $-18808030478336+3274057859072\sqrt{33}   $\\\cline{2-2}
           & $-45012+7836\sqrt{33},~-5198438+904932\sqrt{33}   $\\\cline{2-2} 
           &$-(87-15\sqrt{33}), ~-2,  $\\\cline{2-2}
           &$-2 $\\\hline 
 $-123  $ & $ x^2+1354146840576\cdot 10^3x+148809594175488\cdot 10^6    $\\\cline{2-2}
           & $-677073420288000+105741103104000\sqrt{41}   $\\\cline{2-2}
           & $-960+120\sqrt{41},~-13314+2240\sqrt{41}   $\\\cline{2-2} 
           &$ -24, ~-2(32+5\sqrt{41})^2 $\\\cline{2-2}
           &$-2 $\\\hline 
 $-147  $ & $x^2+34848505552896\cdot 10^3x+11356800389480448\cdot 10^6 $\\\cline{2-2}
           & $-17424252776448000+3802283679744000\sqrt{21}   $\\\cline{2-2}
           & $-2520-240\sqrt{21},~-31724-11418\sqrt{21}   $\\\cline{2-2} 
           &$63+9\sqrt{21}, ~(7-\sqrt{21})((5+\sqrt{21})/2)^8  $\\\cline{2-2}
           &$7-\sqrt{21} $\\\hline 
 $-267  $ & $ x^2+19683091854079488\cdot 10^6x+531429662672621376897024\cdot 10^6 $\\\cline{2-2}
           & $-9841545927039744000000+1043201781864732672000\sqrt{89}   $\\\cline{2-2}
           & $-37500+3180\sqrt{89},~3250002-371000\sqrt{89}   $\\\cline{2-2} 
           &$150, ~2(500+53\sqrt{89})^2  $\\\cline{2-2}
           &$2 $\\\hline 
\multicolumn{2}{l}{}\\
\multicolumn{2}{l}{The results and data for the case $h(R)=3$}\\
\hline           
$-108  $ & $x^3-151013228706\cdot 10^3x^2+224179462188\cdot 10^6x$\\
&$\phantom{aaaaaaaaaaa}-1879994705688\cdot 10^9$\\\cline{2-2}
           & $31710790944000\sqrt[3]{4}+39953093016000\sqrt[3]{2}+50337742902000   $\\\cline{2-2}
           & $105\sqrt[3]{4}-90\sqrt[3]{2}-135,~-738\sqrt[3]{4}+738\sqrt[3]{2}+526   $\\\cline{2-2} 
           &$9-3\sqrt[3]{2}, ~4(1-\sqrt[3]{2})^8(-1+\sqrt[3]{2})  $\\\cline{2-2}
           &$-1+\sqrt[3]{2} $\\\hline
$-243  $ & $ x^3+1855762905734664192\cdot 10^3x^2-3750657365033091072\cdot 10^6x$\\
         &$+3338586724673519616\cdot 10^9    $\\\cline{2-2}
           & $-618587635244888064000-428904711070941184000\sqrt[3]{3}$\\
           &$-297385917043138560000\sqrt[3]{9}   $\\\cline{2-2}
           & $-1560+720\sqrt[3]{9},~32258-11124\sqrt[3]{3}-7704\sqrt[3]{9}   $\\\cline{2-2} 
           &$42-18\sqrt[3]{9}, ~(-2+\sqrt[3]{9})^8(-4+2\sqrt[3]{9})  $\\\cline{2-2}
           &$-4+2\sqrt[3]{9} $\\\hline             
\end{tabular}\vspace{5mm}
\end{center}
 {\bf 3.2.} The case $s=4$\vspace{3mm}\par
 In this case, by Lemma~\ref{lem3} and Proposition~\ref{prop4}, we know that $x_Q$ is one of two $F$-rational solutions of $\Psi_4(x,E)/y=0$ satisfying the condition given in the last part of Lemma~\ref{lem3}. We shall explain the case $d(R)=-32$. 
We take a solution $j=26125000+18473000\sqrt 2$ of $H_{32}(x)=x^2-52250000x+12167000000=0$ and consider an elliptic curve $E$ with $j(E)=j$, defined by an equation:~
$$y^2=x^3+Ax+B~\quad (A =-105-90\sqrt{2},~B=630+518\sqrt{2}).$$
Then $\Psi_4(x,E)/y=0$ has two $F$-rational solutions $x_1=3+5\sqrt 2,x_2=9-\sqrt 2$. Consider a prime number $p=17=3^2+2^2 2$ and a prime ideal $\p=(1-3\sqrt 2)$. Then by counting the number of points of $\overset{\sim}{E_\p}(\mathbb F_p)$, we know $a_\p(E)=-6$. Since $(x_1^3+Ax_1+B/\p)=(-3+3\sqrt 2/\p)=(-2/17)=1=(-1)^{(a_\p(E)/2-1)}$, we see $x_Q=x_1$. Calculating $y_Q^2$, we may obtain $\alpha_E=-3+3\sqrt{2}$.
For the cases $d(R)=-64,-112$, we know the class polynomials are:
\begin{equation*}
\begin{cases}
H_{64}(x)=x^2-82226316240x-7367066619912,\\
H_{112}(x)=x^2-274917323970000x+1337635747140890625. 
\end{cases}
\end{equation*}
In these cases by similar argument we have $(E,\alpha_E)$. We list our results in the next propositon.
\begin{prop}\label{prop6}\phantom{aaa}
\end{prop}
\begin{center}
\begin{tabular}{|l|l|l|l|}\hline
$d(R)$ & $j(E)$ & $A$ & $\alpha_E$\\
       &&$B$&\\\hline
$-32$ & $26125000+18473000\sqrt 2$ & $-105-90\sqrt{2}$& $-3+3\sqrt{2}$ \\
&&$630+518\sqrt{2}$&\\\hline
$-64$ & $41113158120+29071392966\sqrt 2$& $-91-60\sqrt 2 $&$\sqrt 2-1$\\
&&$ 462+308\sqrt 2$&\\\hline
$-112$ & $137458661985000$&$-725-240\sqrt 7$&$1$\\
&$\phantom{aa}+51954490735875\sqrt 7$&$9520+3698\sqrt 7$ &\\\hline
\end{tabular}
\end{center}

\noindent
{\bf 4.}~ The case $mf^2$ is divided by $5$\vspace{3mm}\par
We shall consider the orders $R$ of discriminant $d(R)=-20$,$-35$,$-40$,\newline $-115$, $-235$. These orders $R$ are maximal and of class number $2$. Further for any $R$, we know $F=\mathbb Q(\sqrt{5})$. For a given order $R$, we consider an elliptic curve $E$, defined over $F$, with complex multiplication by $R$. Proposition~\ref{prop3} shows that $\Psi_5(x,E)$ has only one $F$-rational factor $H_{1,E}(x)$ of degree $2$ and for any solution $x_1$ of $H_{1,E}(x)=0$, a point $Q$ of $E$ with $x_Q=x_1$ is a generator of the group $C_5$. Let $L=F(x_Q)$ and $M=F(Q)$. For a prime number $p$ satisfying $p^{\ell_\p}\modx 15$, our problem is rather easy (cf. Proposition~\ref{prop2}).  For a prime number $p$ satisfying $p^{\ell_\p}\modx 45$ and a prime ideal $\p$ dividing $p$, to determine $r_\p$, we must study the ideal class groups of $F$ corresponding to the fields $L$ and $M$. Regarding conductors of $L$ and $M$, we have a following result. In Proposition~\ref{prop7}, we shall use the notation in \S 2.
\begin{prop}\label{prop7} Let $\q$ be a prime ideal of $F$ prime to $(2\sqrt 5)$. Then $e_\q(\f_L)=e_\q(\f_M)$ and $e_\q(\f_M)\le 1$. Further if $E$ has good reduction at $\q$, then $e_\q(\f_M)=0$.
\end{prop} 
\begin{proof} Since $M$ is a cyclic extension of degree $4$ over $F$, we have $e_\q(\f_M)\le 1$  (cf. Serre [9]). The other assertions are deduced from Lemma~\ref{lem4}. 
\end{proof}
As for the prime ideal $(\sqrt 5)$, we have $e_{(\sqrt 5)}(\f_L)\le e_{(\sqrt 5)}(\f_M)\le 1$. 
  Proposition 7 shows, to avoid tedious computation in determining class groups, it is necessary to choose an elliptic curve $E$ such that the number of prime factors of its discriminant is as small as possible. 
 We shall explain the case $d(R)= -235$, because the other cases can be deduced from similar but much easier argument.  We have
 $$H_{235}(x)=x^2+82317741944942592\cdot 10^4x+11946621170462723407872\cdot 10^3.$$
 First we take a solution $$j=-411588709724712960000-184068066743177379840\sqrt{5}$$ of $H_{235}(x)=0$.

We consider an elliptic curve $E$ defined by an equation:~$y^2=x^3+(-15510+2068\sqrt{5})x+(3200841-649446\sqrt{5})/4$. The discriminant of $E$ is $47^32^{-4}e^{-42}$, $j(E)=j$ and   
 $$H_{1,E}(x)=10x^2+(3525-2115\sqrt{5})x+624160-262918\sqrt{5},$$
 where $e=(1+\sqrt 5)/2$. By solving the equation $H_{1,E}(x)=0$, we obtain a generator $Q$ of $C_5$ given by
 $$ x_Q=3e^{-1} t+47 e^{-10}/2,~ y_Q=(2\sqrt 5e^{10})^{-1}\pi$$,
  where $t=\sqrt{47\sqrt{5}e^{-1}}$ and $\pi=\sqrt{47e^{-1}(2115-(211+23\sqrt 5)t)}$.
In particular we have $L=F(t)$ and $M=L(\pi)$. Next we shall determine conductors and ideal class groups of $L$ and $M$. Since the maximal order of $L$ has a basis $\{1,(1+e^{-1}t)/2\}$ over the maximal order of $F$, the discrimant of $L$ over $F$ is $(e^{-1}t)^2$. This shows that $\f_L=(47\sqrt{5})$. Since $M$ is real and a conjugate field of $M$ over $\mathbb Q$ is imaginary, Proposition~\ref{prop7} shows $\f_M=(2^k\cdot47\sqrt{5})\infty_2$, where $k\in\mathbb Z, 0\le k\le 2$ and $\infty_2$ is the infinite place of $F$ corresponding to the conjugate embedding of $F$ to $\overline{\mathbb Q}$. We have only to determine the $2$-exponent $k$. See \S 3 of Hiramatsu and Ishii [4] for a method to calculate the 2-exponent of conductors. For a moment, we assume $M$ is defined modulo $(4\cdot47\sqrt{5})\infty_2$.  Let $\P$ be the ray class group of $F$ modulo $(4\cdot 47\sqrt{5})\infty_2$. Denote by $\P_L$ and $\P_M$ the subgroups of $\P$ corresponding to $L$ and $M$ respectively. Consider the ideal classes $\g,\k$ and $\l$ of $\P$ represented by the principal ideals $((1+3\sqrt{5})/2),(46+47\sqrt{5})$ and $(471)$ respectively. Then we see $\g$ is of order $276$ and both $\k$ and $\l$ are of order $2$. Further we have:
$$\P= <\g>\times <\k>\times <\l>\quad \text{(a direct product)}.$$
 Let $\P_1$ be the ray class group modulo $(47\sqrt{5})$ and $\theta$ the canonical morphism of $\P$ to $\P_1$. Then $\P_1$ is a cyclic group generated by $\theta(\g)$ of order $138$ and $\rit{Ker}(\theta)=<g^{138},\k,\l>$. Since $f_L=(47\sqrt 5)$, $\P_L\supset \rit{Ker}(\theta)$. This shows that $\P_L=<\g^2,\k,\l>$. Next we shall determine $\P_M$. Let $\xi$ be the homomorphism of $\P$ to itself defined by $\xi (\mathfrak a)=\mathfrak a^{69}$. Then $\xi$ induces an isomorphism of $\P/\P_M$ to $\xi(\P)/\xi(\P_M)$. Consider the prime numbers $q_1=251=4^2+235,q_2=431=14^2+235$ and $q_3=239=2^2+235$ and prime ideals $\q_1=(16+\sqrt{5}),\q_2=((43+5\sqrt{5})/2)$ and $\q_3=((31+\sqrt{5})/2)$ of $F$ dividing $q_1,q_2$ and $q_3$ respectively. In the following, for a prime ideal $\q$ of $F$, we denote by $C(\q)$ the class of $\P$ represented by $\q^{69}$. Then we know  $C(\q_1),C(\q_2)$ and $C(\q_3)$ belong to $\k\l,~\k$ and $\xi(\g)\k$ respectively. By counting the number of rational points of the reduced elliptic curve of $E$ modulo $\q_i$, we have $a_{\q_1}(E)=-8,a_{\q_2}(E)=-28$ and $a_{\q_3}(E)=-4$. Therefore, by Theorem~\ref{th1}, we know $\k,\l\in \P_M$ and the class $\xi(\g)\k$ corresponds to the isomorphism $\lambda$ such that $Q^\lambda =[3](Q)$. Since $\P_M$ is a subgroup of $\P_L$ of index $2$, we conclude that $\P_M=<\g^4,\k,\l>$. In particular, $\P_M$ does not contain the kernel $<\g^{138}\k,\l>$ of the canonical morphism of  $\P$ to the ray class group modulo $(2\cdot 47\sqrt{5})\infty_2$. Therefore $\f_M=(4\cdot 47\sqrt{5})\infty_2$. We know the class $\m=\xi(\g)$ is represented by the ideal $(743+756\sqrt{5})$. Consequently we have
\begin{thm} Let $\k,\l$ and $\m$ be the classes of $\P$ represented by the ideals $(46+47\sqrt{5})$, $(471)$ and $(743+756\sqrt{5})$ respectively. Put $\SS=<\m,\k,\l>$ and $\D=<\k,\l>$. Let $p$ be an odd prime number and $\p$ a prime ideal of $F$ dividing $p$ and assume that they satisfy {\rm (1)}. Furthermore, let $u_p$ and $v_p$ be the positive integers such that $4p^{\ell_\p}=u_p^2+235v_p^2$ and $(u_p,p)=1$. If the class $C(\p)$ of $\p^{69}$ belongs to $\m^i\D~(0\le i\le 3)$, and $\epsilon_\p\in\{\pm 1\}$ is chosen such that $\epsilon_\p u_p \modx {2\cdot 3^i} 5$, then we have  $a_\p(E)=\epsilon_\p u_p$.
\end{thm}
\noindent
\begin{note}
$C(\p)\in \D\cup\m^2\D$ if and only if $p^{\ell_\p}\modx 15$.
\end{note}
\begin{exam} 
\begin{enumerate}
\itemx i Let $p=239=2^2+235$ and $\p=((31+\sqrt 5)/2)$. Then $C(\p)=\m\k\in \m D$ and $a_\p(E)=-4$.
\itemx {ii} Let $p=241=(27^2+235)/4$ and $\p=((33+5\sqrt 5)/2)$. Then $C(\p)=\l\in D$ and $a_\p(E)=27$.
\itemx {iii}Let $p=719=22^2+235$ and $\p=((59+11\sqrt 5)/2)$. Then $C(\p)=\m^3\k\l \in \m^3D$ and $a_\p(E)=44$.
\end{enumerate}
\end{exam}
We shall give the data and results for other cases. In the below, put $t_m=\sqrt{m/(\sqrt{5}e)}$ and denote by $\P$ the ray class group of conductor $\f_M$ of $F$. Further, we denote by $p$ and $\p$ an odd prime number and a prime ideal of $F$ dividing $p$ such that they satisfy the condition (1) for the given elliptic curve $E$.\vspace{3mm}\newline
(I)~The case $m=5,~d(R)=-20$

\begin{equation*}
\begin{cases}
H_{20}(x)=x^2-1264000x-681472000,\quad j(E)=632000+282880\sqrt{5},\\
A=-50/3-5\sqrt{5},\quad B=100/3+280\sqrt{5}/27, \\
x_Q=5e^2/6+t_m, ~y_Q=(\sqrt{5})(e+t_m)\sqrt{1+t_m^{-1}}, \\
 L=F(t_m),M=L\left (\sqrt{1+t_m^{-1}}\right ),~\f_L=(4\sqrt{5}),~\f_M=(8\sqrt{5}),\\
\P= <\g_1>\times <\g_2>,~\P_L=<\g_1^2,\g_2>,~ \P_M=<\g_1^2\g_2>,
\end{cases}
\end{equation*}
where $\g_1$ and $\g_2$ are the classes  of order of $4$ and $2$ represented by the ideals $((21+\sqrt{5})/2)$ and $(11+2\sqrt{5})$. 
\begin{prop} Let $u_p$ and $v_p$ be the positive integers such that $p^{\ell_\p}=u_p^2+5v_p^2,~(u_p,p)=1$. Choose $\epsilon_\p\in\{\pm 1\}$ such that $\epsilon_\p u_p\modx {2^i}5$ if the class of $\p$ belongs to $\g^i_1\P_M ~(0\le i\le 3)$. Then we have $a_\p(E)=2\epsilon_\p u_p$.\end{prop}
 Choosing a suitable generator of $\p$, $p$ is written in a form $p=a^2-5b^2$, where  $a$ and $b$ are integers satisfying the condition:
\begin{align*}
a&\equiv\begin{cases}17 \mod{20}~&\text{ if }p\modx 45 \\
   \phantom{a}1\mod 20 ~&\text{ if }p\modx 15 ,\end{cases} \\
b&\equiv\begin{cases} 2\mod 4~&\text{ if } p\modx 58\\
0 \mod 4~&\text{ if }p\modx 18.\end{cases}
\end{align*} 
For $i=1,2$, let $p_i=a_i^2-5b_i^2$ be the prime numbers represented as above. If $p_1\modx {p_2}{40}$, then prime ideals $(a_1+b_1\sqrt{5})$ and $(a_2+b_2\sqrt{5})$ belong to the same class of $\P$ if and only if $a_1-a_2+5(b_1-b_2)\modx 0{40}$. Let $\T=<\g_1^2>$. Then we see if $p\modx{1~(\text{resp.}~9,21,29)}{40}$, then the class $C(\p)$ of the prime ideal $\p=(a+b\sqrt 5)$ belongs to $\T,(\text{resp.}~\g_2\g_1\T,\g_2\T,\g_1\T)$ and furthermore $C(\p)\in P_M$ if and only if $a+5b\modx{1 (\text{resp.}~-3,11,7)}{40}$. Therefore we have
\begin{prop} Let $p=u^2+5v^2=a^2-5b^2$, where $u$ and $b$ are positive integers and $a$ and $b$ are integers satisfying the above condition. Then if we choose $\epsilon_\p\in\{\pm 1\} $such that
\begin{equation*}
\epsilon_\p u\equiv \begin{cases}
(-1)^{(a+5b-1)/20}a\mod 5\quad &\text{if }p\modx {1}{40},\\
(-1)^{(a+5b+3)/20}a\mod 5\quad &\text{if }p\modx {9}{40},\\
(-1)^{(a+5b-11)/20}a\mod 5\quad &\text{if }p\modx {21}{40},\\
(-1)^{(a+5b-7)/20}a\mod 5\quad &\text{if }p\modx {29}{40},
\end{cases}
\end{equation*}
then we have $a_\p(E)=2\epsilon_\p u$.
\end{prop}
\noindent
(II)~The case $m=35,d(R)=-35$

\begin{equation*}
\begin{cases}
H_{35}(x)=x^2+117964800x-134217728000,\\
j(E)=-58982400-26378240\sqrt 5,\\
A=-70\sqrt 5/3,~B=(13475+980\sqrt 5)/108,\\
x_Q=(35e-3t_m)/6e^2, ~y_Q=(t^2_m/2e^2)\sqrt{\sqrt 5-(9+\sqrt 5)(et_m)^{-1}},\\
L=F(t_m),M=L\left (\sqrt{\sqrt 5 - (9+\sqrt 5)(et_m)^{-1}}\right ),\\ 
\f_L=(7\sqrt{5}),~\f_M=(14\sqrt{5}\infty_2),\\
\P =<\h>,~\P_L=<\h^2>,~ \P_M=<\h^4>,
\end{cases}
\end{equation*}
 where $\h$ is the class of order $12$ represented by the ideal $(6+\sqrt{5})$ . 
\begin{prop}  Let $u_p$ and $v_p$ be the positive integers such that $4p^{\ell_\p}=u_p^2+35v_p^2,~(u_p,p)=1$.   Choose $\epsilon_\p\in\{\pm 1\}$ such that $\epsilon_\p u_p \modx {2\cdot 3^i} 5$ if the class of $\p$ belongs to $\h^i\P_M ~(0\le i\le 3)$. Then we have $a_\p(E)=\epsilon_\p u_p$.
\end{prop}
\noindent
(III)~The case $m=10,~d(R)=-40$

\begin{equation*}
\begin{cases}
H_{40}(x)=x^2-425692800x+9103145472000,\\
j(E)=212846400+95178240\sqrt{5},\\
A=-125+15\sqrt 5, B=-200+240\sqrt{5},\\
x_Q=(10e+t_m)/e^2,~y_Q=2e^{-2}t_m\sqrt{15e^{-1}+(40-11\sqrt 5)t^{-1}_m},\\
L=F(t_m),M=L\left(\sqrt{15e^{-1}+(40-11\sqrt 5)t^{-1}_m}\right),\\ 
\f_L=(8\sqrt{5}),~\f_M=(16\sqrt{5)}\infty_2,\\
\P =<\g>\times <\h>\times <\l>,~\P_L=<\g^2,\h,\l>,~ \P_M=<\h,\l>,
\end{cases}
\end{equation*}
 where $\g, \h$ and $\l$ are the classes of order $4,4$ and $2$ represented by the ideals $(6+\sqrt 5)$,$((53+3\sqrt{5})/2)$ and $((37+7\sqrt{5})/2)$ respectively. 
\begin{prop} Let $u_p$ and $v_p$ be the positive integers such that $p^{\ell_\p}=u_p^2+10v_p^2,~(u_p,p)=1$. Choose $\epsilon_\p\in\{\pm 1\}$ such that $\epsilon_\p u_p \modx {2^i}5$ if the class of $\p$ belongs to $\g^{i}\P_M~ (0\le i\le 3)$. Then we have $a_\p(E)=2\epsilon_\p u_p$.
\end{prop}
\noindent
(IV)~The case $m=115,~d(R)=-115$
\begin{equation*}
\begin{cases}
H_{115}(x)=x^2+427864611225600x+130231327260672000,\\
j(E)=-213932305612800+95673435586560\sqrt{5},\\
A=-345-23\sqrt{5}, B=-(19573+5290\sqrt{5})/4,\\
x_Q=e^3\sqrt{5}t_m(\sqrt{5}t_m+e^4)/10,\\
y_Q=(e^9t_m^2/10)\sqrt{15e^{-1}+(-85+61\sqrt{5})t^{-1}_m},\\ 
L=F(t_m),M=L\left(\sqrt{15e^{-1}+(-85+61\sqrt{5})t^{-1}_m}\right),\\
\f_L= (23\sqrt{5}),~\f_M=(92\sqrt{5})\infty_2,\\
\P =<\f_1>\times <\f_2>\times <\f_3>,~\P_L=<\f_1^2,\f_2,\f_3>,~ \P_M=<\f_1^4,\f_2,\f_3>,  \end{cases}
\end{equation*}
where $\f_1,\f_2$ and $\f_3$ are the classes represented by the ideals $((1+3\sqrt{5})/2),(24+23\sqrt{5})$ and $(91)$ and the order of $\f_1,\f_2$ and $\f_3$ are $132,2$ and $2$ respectively. Since the map $\xi_1:\mathfrak a \to \mathfrak a^{33}$ of $\P$ to itself induces an isomorphism of $\P/\P_M$ to $\xi_1(\P)/\xi_1(\P_M)$
and $\f_0=\f_1^{33}$ is represented by the ideal $(423+372\sqrt{5})$, we have 
\begin{prop}  Let $\SS=<\f_0,\f_2,\f_3>$ and $\D=<\f_2,\f_3>$. Let $u_p$ and $v_p$ be the positive integers such that $4p^{\ell_\p}=u_p^2+115v_p^2,~(u_p,p)=1$.   Choose $\epsilon_\p\in\{\pm 1\}$ such that $\epsilon_\p u_p \modx {2\cdot 3^i} 5$  if the class of $\p^{33}$ belongs to $\f_0^i\D~(0\le i\le 3)$. Then we have $a_\p(E)=\epsilon_\p u_p$.
\end{prop}
\vspace{3mm}\noindent
{\it{\bf Acknowledgement.} The author would like to express his hearty gratitude to Professor M. Kaneko for offering the table of class equations and to Tomoko Ibata for writing and running the computer program required for this work.}

Noburo Ishii \vspace{3mm}\newline
 Department of Mathematics and Information Science,\newline 
Osaka Prefecture University, \newline 
1-1 Gakuen-cho, Sakai, Osaka 599-8531 \newline 
Japan \newline
e-mail:ishii@mi.cias.osakafu-u.ac.jp  
\end{document}